\documentclass{amsart}
\usepackage{longtable}
\usepackage{hyperref}
\usepackage{amsrefs}
\usepackage{comment}
\usepackage{float}
\usepackage[dvipsnames]{xcolor}
\usepackage{tikz}
\usepackage{multirow}
\usepackage{booktabs}
\usepackage{caption}

\theoremstyle{definition}
\newtheorem{theorem}{Theorem}

\newcommand{\tri}{\mathbin{\overline{\nabla}}}

\newcommand{\fub}{ }

\newcommand{\m}{{\mathbf m}}
\newcommand{\C}{{\mathbf c}}
\newcommand{\T}{{\mathbf t}}

\newcommand{\Sb}{{\mathbf S}}

\newcommand{\Tb}{{\mathbf T}}
\newcommand{\Mb}{{\mathbf M}}
\newcommand{\Ms}{{\mathcal M}}

\newcommand{\Sm}{\mathcal{S}} 

\newcommand{\mC}{C_{\m}}
\newcommand{\mV}{V_{\m}}
\newcommand{\mE}{E_{\m}}
\newcommand{\mF}{F_{\m}}

\newcommand{\BigComment}[1]{}

\newcommand{\betaALL}{\beta}
\newcommand{\betaSOME}{\beta_\textrm{L}}
\newcommand{\kmod}{~\mathrm{ mod }~}
\newcommand{\kmodf}{\kmod f^{d+1}}
\newcommand{\bbb}{\!\!\!\!\!\!\!\!}

\newcommand{\n}{{\bf n}}

\providecommand{\U}[1]{\protect\rule{.1in}{.1in}}

\begin{document}
\title[ Finite Interpretation of the Hyper-Catalan Series Zero Powers]{Finite Interpretation of the Hyper-Catalan Series Zero and its Powers}
\author{Dean Rubine and Pratham Mukewar}

\date{\today}
\begin{abstract}
In 2025,
Wildberger and Rubine showed the formal series zero of the univariate geometric polynomial 
is $\Sb$, the generating series for the hyper-Catalan numbers $\mC$,
which count the number of roofed subdivided polygons (subdigons) of type $\m$.
We show that we can interpret this result as a finite identity at each level, where a level is a truncation of $\Sb$ to a given maximum number of vertices, edges, or faces (bounded by degree) of the associated subdigon types.
We then explore powers $\Sb^r$, recounting Raney's and our own combinatorial derivations of its coefficients.

\end{abstract}
\maketitle
\section{Introduction}
Wildberger and Rubine~\cite{Wildberger2025} call
\begin{equation} \label{eqn:g}
g(\alpha) = 1- \alpha + t_2 \alpha^2 + t_3 \alpha^3 + t_4 \alpha^4 + \ldots \end{equation}
the \textbf{general geometric polynomial},
which means general except for the constant being $1$ and the linear coefficient $-1$.
In contrast, their formal series zero of $g(\alpha)$, 
\begin{equation}
\Sb 
\equiv \sum_{m_2, m_3, \ldots \ge 0} C[m_2, m_3, \ldots]  t_2^{m_2} t_3^{m_3} t_4^{m_4} \cdots
\equiv \sum_{\m \ge 0} \mC \T^{\m},   
\end{equation}
is an ongoing multivariate series in a potentially unbounded number of variables $t_2, t_3, \ldots$, so more difficult to grasp.

Wildberger's result extends to univariate power series, but
here we are interested in the case where $g(\alpha)$ is a polynomial, 
a finite object whose zeros we seek.
That way, for $\alpha_L$, some finite truncation of $\Sb$,
the evaluation $g(\alpha_L)$ is a finite computation we can completely perform.

\begin{figure}
    \centering
    \includegraphics[width=1\linewidth]{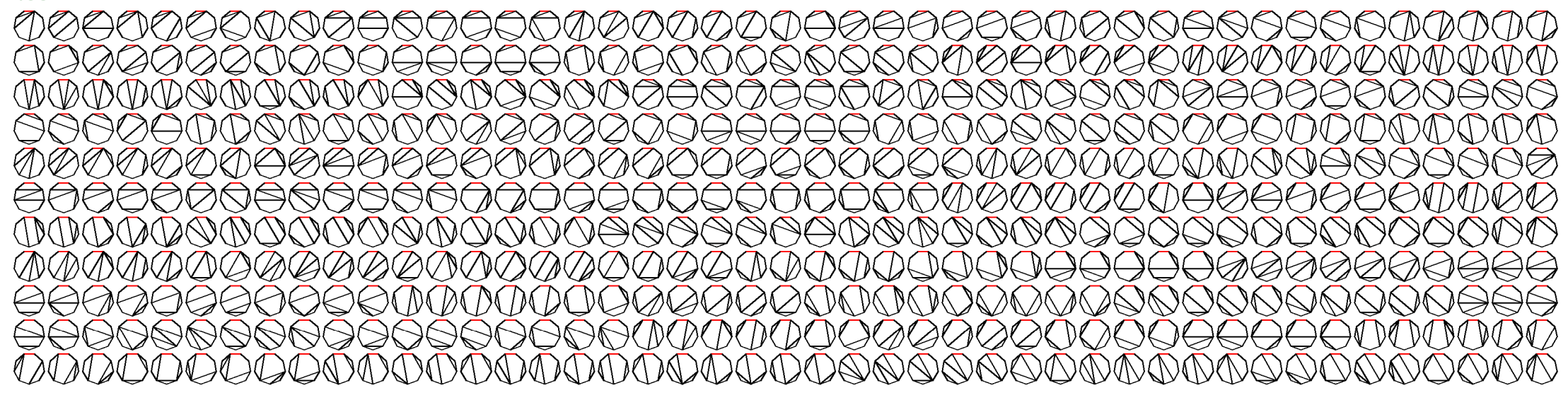}
    \caption{C[2,1,1]=495 subdigons with two triangles, one quadrilateral and one pentagon}
    \label{fig:C211}
\end{figure}

$\mC$ is the array of \textbf{hyper-Catalan} numbers,
which count the number of subdigons with $m_2$ triangles, $m_3$ quadrilaterals, etc., i.e. the number of subdigons of \textbf{type} $\m=[m_2, m_3, \ldots]$ (see figure \ref{fig:C211} and section \ref{sec:subdigons}).
Its closed form~\cites{Schuetz2016,Rhoades2011,Erdelyi1940} hints how we might organize layers of $\Sb$:
\begin{align} \label{eqn:CVMF}
\mC = \dfrac{( \mE-1) !}{ (\mV-1)! \, \m!},  \qquad \textrm{ where } 
\m! \equiv  m_2! \, m_3! \, m_4! \, \cdots \textrm{ and } 
\\ \nonumber
\quad
\mV
= 2+\sum_{k \ge 2}^{\fub} (k-1)m_k,
\quad 
\mE 
= 1 + \sum_{k \ge 2}^{\fub} k \, m_k ,
\quad
\mF 
= \sum_{k \ge 2}^{\fub} m_k
\end{align}

$\mV$, $\mE$ and $\mF$ respectively count of the number of vertices, edges and faces of a subdigon of type $\m$.
$\mF$ doesn't appear explicitly in the formula for $\mC$, but is implicit in the factor of $\m!$.
The formulas satisfy the Euler Characteristic relation $\mV - \mE + \mF =1$.

A subdigon $s$ of type $\m$ has an \textbf{accounting monomial} $\psi(s) \equiv \T^\m$, where the type vector $\m$ is recorded as the powers of the $t_k$.  
\begin{equation}
    \Psi(M) \equiv \sum_{s \in M} \psi(s)
\end{equation}is the polynomial associated with a multiset of subdigons $M$, formed as the sum of the accounting monomials of all the subdigons in $M$.
For the multiset $\Sm_{[2,1,1]}$ pictured in figure \ref{fig:C211}, all the accounting monomials are identical, so we have $\Psi(\Sm_{[2,1,1]})=495 t_2^2 t_3 t_4$.

\section{Layering}
 \label{sec:layering}
 
We add layering variables for vertices, edges and faces to our geometric polynomial equation and solution, similarly encoding the counts in the exponents. Since the null subdigon $|$ (two vertices, no faces, one edge), has accounting monomial $1$, it is best to make the power of $v$ \textit{two less} than the number of vertices, and the power of $e$ \textit{one less} than the number of edges. 

\begin{theorem}[Series Zero of the Layered Polynomial] \label{thm:layered}
Given
$$  h(\alpha) = 1 -\alpha + \sum_{n \ge 2}^{\fub} t_n  v^{n-1} e^{n} f^1 \alpha^n , \qquad 
\betaALL = \sum_{\m \ge 0} \ \mC v^{\mV-2}e^{\mE-1}f^{\mF}   \T^{\m} 
 $$
we have  $ h(\betaALL) = 0 $.
\end{theorem}

\begin{proof}
The coefficient on $\alpha^n$ in $h(\alpha)$ is $c_n = t_n  v^{n-1} e^{n} f^1$ 
for $n \ge 2$.
Note:
\begin{align}
 \C^{\m} 
& = \prod_{n \ge 2} c_n^{m_n}
=\prod_{n \ge 2} \left(t_n  v^{n-1} e^{n} f^1 \right) ^{m_n}
\\ &
= v^{   \left( \sum_{n \ge 2} (n-1)m_n \right)} e^{  \left( \sum_{n \ge 2} n m_n \right) } f^{  \left( \sum_{n \ge 2} m_n \right)}  \prod_{n \ge 2} t_n ^{m_n}
=  v^{ \mV-2 } e^{\mE-1} f^{\mF}  \T^{\m}
\nonumber 
\end{align}
so we have $  \beta =  \sum_{\m \ge 0} \ \mC \,  \C^{\m} $.  Thus
$ h(\betaALL)  =  1 -\betaALL + \sum_{n \ge 2}^{\fub} c_n \betaALL^n 
   = 0
$
by Wildberger and Rubine's soft geometric polynomial formula.
\end{proof}
We're projecting the subdigon multiset algebra, at the level of types.  Adding the layering variables doesn't change the miraculous cancellation of the geometric polynomial solution.

$\betaALL$ is ongoing, a multivariate sum.
We're interested in layering the sum, with each layer having a finite number of terms.
We'll use mod in the sense of Knuth~\cite{Graham1989}, as the remainder operator,
as another way to write layerings of the general series $\betaALL$.  
We define:
\begin{align}
\beta_{d.} & \equiv  \sum_{\mV-2 \le d}  \! \mC \, v^{\mV-2}e^{\mE-1}f^{\mF}   \T^{\m}  = \betaALL \kmod v^{d+1}
\\ 
\beta_{d\_\!\_} & \equiv  \sum_{\mE-1 \le d} \mC  \, v^{\mV-2}e^{\mE-1}f^{\mF}   \T^{\m} = \betaALL \kmod e^{d+1}
\\
\beta_{d\Delta} & \equiv \sum_{ \mF \le d} \mC \, v^{\mV-2}e^{\mE-1}f^{\mF}   \T^{\m} = \betaALL \kmod f^{d+1} 
\\ 
\beta_{d\Delta q} & \equiv \sum_{\substack{\m=[m_2,...,m_q] \ge 0 \\ \mF \le d}}  \mC \, v^{\mV-2}e^{\mE-1}f^{\mF}   \T^{\m} 
\end{align}

The dot, underscore and triangle in the $\beta$ subscripts refer to vertex, edge and face layers respectively.  
$\beta_{d \Delta}$ is an ongoing power series; there is no limit to the number of subdigons with even just a single face.
$\beta_{d\Delta q}$ is finite, a multivariate polynomial, restricted to no more than $(q+1)$-gons, suitable for solving equations of degree $q$ or less.
$\beta_{d.}$ and $\beta_{d\_\!\_}$ are finite as well.
If $\betaSOME$ is one of these finite sums, and $h(\alpha)$ a polynomial, then $h(\betaSOME)$ is a finite computation that we wish to understand.

We need a theorem from the modular arithmetic of polynomials, that evaluating a polynomial on a remainder is the same as evaluating it on the original value,
as long we only compare remainders at the end, which we notate:
\begin{equation}
h(\beta \kmod  z^{d+1} ) \kmod z^{d+1}  = h(\beta) \kmod  z^{d+1} 
\end{equation}
For bounded face layers,
our equation $h_q(\alpha)$ is a $q$ degree polynomial. Let:
\begin{equation}
h_q(\alpha) = 1 -\alpha + \sum_{n= 2}^{q} t_n  v^{n-1} e^{n} f^1 \alpha^n 
, \quad 
\beta_q = \bbb \sum_{\m=[m_2,...,m_q] \ge 0}  \bbb \mC \, \T^{\m}  v^{\mV-2} e^{\mE-1} \,f^{\mF} 
\end{equation}

\begin{theorem} (Vertex, Edge and Face Layering Theorems)  Let $h(\alpha), h_q(\alpha),$ $\betaALL, \beta_{d.},$ $ \beta_{d\_}, $ $ \beta_{d\Delta}$ and $\beta_{d\Delta q}$ be defined as above.  Then:
\begin{align*}
h(\beta_{d.}) \kmod v^{d+1} &= 0, \qquad  h(\beta_{d\_\!\_}) \kmod e^{d+1} = 0,
\\  h(\beta_{d\Delta} ) \kmod f^{d+1} &= 0, \qquad h_q(\beta_{d\Delta q} )\kmodf =0  .
\end{align*}
\end{theorem}
\begin{proof}
The straightforward proofs~\cite{Mukewar2025} use the previous idea about remainders,
as well as Theorem \ref{thm:layered}, that $h(\betaALL)=0$.
\begin{align*}
 h(\beta_{d.})\kmod v^{d+1} &= h(\betaALL \kmod v^{d+1}) \kmod v^{d+1}  = h(\betaALL) \kmod v^{d+1} =  0
\\ h(\beta_{d\_\!\_})\kmod e^{d+1} &= h(\betaALL \kmod e^{d+1}) \kmod e^{d+1}  = h(\betaALL) \kmod e^{d+1} =  0
\\ h(\beta_{d\Delta}) \kmod f^{d+1}&= h(\betaALL \kmod f^{d+1}) \kmod f^{d+1}  = h(\betaALL) \kmod f^{d+1} =  0
\end{align*}
Constraining $\beta_q$ to only have nonzero $m_2$ through $m_q$ doesn't hinder it from being the zero of $q$ degree polynomial $h_q(\alpha)$.
From our definition, we see $ \displaystyle \beta_{d\Delta q}  = \beta_q  \kmod f^{d+1}.$
It follows:
$$
h_q(\beta_{d\Delta q} )\kmodf 
= h_q( \beta_q \kmodf) \kmodf
= h_q( \beta_q ) \kmodf
=0
$$
\end{proof}

The solution to the general geometric polynomial 
is an ongoing, multivariate power series.
In this section we've given a finite interpretation of the layered, ongoing solution $\betaALL$.
For every level $d$, we can take vertex layers of $\betaALL$ up to $d$, and we find that $h$ applied to the layers gives zero, to vertex level $d$.  
We've shown similar results hold for edge layers and bounded face layers.
We can relate the results back to series zeros of $g(\alpha)$ by setting the layering varibles to unity.

\section{Explicit Examples of Layers of {\bf S}}

Let's demonstrate some finite identities by looking at some layerings, starting with vertices. 

\subsection{Vertex Layers}
We'll set $e=f=1$ in $\beta_{5.}$ and $h(\beta_{5.})$ to reduce the clutter.
That's $  h(\alpha) = 1 -\alpha + \sum_{n \ge 2}^{\fub} t_n  v^{n-1}  \alpha^n $. 

We'll start by showing $h(\beta_{5.}) \mod v^6 = 0$.  In fact, for each vertex level $n=0$ to $n=5$ we'll show that 
$ 
[v^n](\beta_{5.}-1) 
=  [v^n]  \sum_{k \ge 2}^{\fub} t_k  v^{k-1}  \beta_{5.}^k
. $ 
 \qquad

Let's save table space on the first few entries:  
\begin{align*}
    &[v^0] (\beta_{5.}-1 ) = [v^0] t_2  v^{1}  \beta_{5.}^2 =0
, \qquad [v^1] (\beta_{5.}-1 ) = [v^1] t_2  v^{1}  \beta_{5.}^2 =t_2
\\& [v^2]  (\beta_{5.}-1 ) =2 t_{2}^{2} + t_{3}
, \quad
[v^2] t_2  v^{1}  \beta_{5.}^2= 2 t_{2}^{2} 
, \quad 
[v^2] t_3  v^{2}  \beta_{5.}^3=  t_{3}.
\end{align*}
\begin{center}
\captionof{table}{Vertex Layering}\label{tab:vertex} \addtocounter{table}{-1}
\begin{longtable}{ | p{2cm} | p{9cm} |} 
\hline
 $\!\![v^{V-2}]$  & $\beta_{5.}-1  = t_2 v \beta_{5.}^2 + t_3 v^2 \beta_{5.}^3 + t_4 v^3  \beta_{5.}^3 + t_5 v^4  \beta_{5.}^5
 + \ldots $
\\
\hline
    \ $ [v^3] t_2  v^{1}  \beta_{5.}^2 $  & $ 5 t_{2}^{3} + 2 t_{2} t_{3}$
 \\ \ $ [v^3] t_3  v^{2}  \beta_{5.}^3 $ & $  \qquad \ \  3 t_{2} t_{3} $
 \\ \ $ [v^3] t_4  v^{3}  \beta_{5.}^4 $ & $  \qquad \qquad \qquad   t_{4}$
 \\   $[v^3] (\beta_{5.}-1 ) $ &  $5 t_{2}^{3} + 5 t_{2} t_{3} + t_{4}$
\\ \hline 
$ \ [v^4] t_2  v^{1}  \beta_{5.}^2 $ & $14 t_{2}^{4} + 12 t_{2}^{2} t_{3} + 2 t_{2} t_{4}   $
\\ $ \ [v^4] t_3  v^{2}  \beta_{5.}^3 $ & $  \qquad \quad \ \   9 t_{2}^{2} t_{3} \qquad \quad\   + 3 t_{3}^{2}$
\\ $\  [v^4] t_4  v^{3}  \beta_{5.}^4 $ & $  \qquad\qquad\qquad \ \ \  4 t_{2} t_{4}$
\\ $ \ [v^4] t_5  v^{4}  \beta_{5.}^5 $ & $  \hspace*{12.5em} \   t_{5}$
\\ $[v^4]  (\beta_{5.}-1 )$ & $14 t_{2}^{4} + 21 t_{2}^{2} t_{3} + 6 t_{2} t_{4} + 3 t_{3}^{2} + t_{5}$
\\ \hline 
 $ \ [v^5] t_2  v^{1}  \beta_{5.}^2 $ & $ 42 t_{2}^{5} + 56 t_{2}^{3} t_{3} + 14 t_{2}^{2} t_{4} + \   7 t_{2} t_{3}^{2} + 2 t_{2} t_{5} $
\\ $ \ [v^5] t_3  v^{2}  \beta_{5.}^3 $ & $  \qquad \quad   28 t_{2}^{3} t_{3}
\hspace{3.8em} + 21 t_{2} t_{3}^{2}
\hspace{3.7em} + 3 t_{3} t_{4}  $
\\ \ $ [v^5] t_4  v^{3}  \beta_{5.}^4 $ & $  \qquad\qquad\qquad \ \ \ 14 t_{2}^{2} t_{4} 
\hspace{7.4em} + 4 t_{3} t_{4} $
\\ $ \ [v^5] t_5  v^{4}  \beta_{5.}^5 $ & $  \hspace{14.5em} 5 t_{2} t_{5} $
\\ $ \ [v^5] t_6  v^{5}  \beta_{5.}^5 $ & $  \hspace{21em} \  t_{6}$
\\ $[v^5]  (\beta_{5.}-1 )$ & $42 t_{2}^{5} + 84 t_{2}^{3} t_{3} + 28 t_{2}^{2} t_{4} + 28 t_{2} t_{3}^{2} + 7 t_{2} t_{5} + 7 t_{3} t_{4} + t_{6}$
\\
\hline
\end{longtable}
\end{center}
We've shown $[v^k] (1 - \beta_{5.} + \sum_{k \ge 2} t_k  v^{k-1}  \beta_{5.}^k) = 0$ for $k=0,1,2,3,4,5$ so indeed $h(\beta_{5.}) \mod v^6 = 0$.
\subsection{Edge Layers}
Similarly to vertices, we set $v=f=1$ in $\beta_{8\_}$ and $h(\beta_{8\_})$.
So, $  h(\alpha) = 1 -\alpha + \sum_{n \ge 2}^{\fub} t_n  e^n  \alpha^n $. 

We will show that $h(\beta_{8\_}) \mod e^9 = 0$, and for each edge level $n=0$ to $n=8$, we'll show that 
$[e^n](\beta_{8\_}-1) 
=  [e^n]  \sum_{k \ge 2}^{\fub} t_k  e^k  \beta_{8\_}^k$. The first few entries are trivial:
\begin{align*}
    &[e^0] (\beta_{8\_}-1 ) = [e^0] t_2  e^2  \beta_{8\_}^2 =0
, \qquad [e^1] (\beta_{8\_}-1 ) = [e^1] t_2  e^2  \beta_{8\_}^2 =0
, \\& [e^2] (\beta_{8\_}-1 ) = [e^2] t_2  e^2  \beta_{8\_}^2 =t_2, \qquad [e^3] (\beta_{8\_}-1 ) = [e^3] t_3  e^3  \beta_{8\_}^3 =t_3
\end{align*}
\begin{center}
\captionof{table}{Edge Layering}\label{tab:edge} \addtocounter{table}{-1}
\begin{longtable}{ | p{2cm} | p{7.4cm} |} 
\hline
 $\!\![e^{E-1}]$  & $\beta_{8\_}-1  = t_2 e^2 \beta_{8\_}^2 + t_3 e^3 \beta_{8\_}^3 + t_4 e^4  \beta_{8\_}^3 + t_5 e^5  \beta_{8\_}^8
 + \ldots $
\\
\hline
    \ $ [e^4] t_2  e^2  \beta_{8\_}^2 $  & $2t_2^2$
 \\ \ $ [e^4] t_4  e^4  \beta_{8\_}^4 $ & \ \  \qquad $t_{4}$
 \\   $[e^4] (\beta_{8\_}-1 ) $ &  $2t_2^2+t_4$
\\ \hline 
$ \ [e^5] t_2  e^2  \beta_{8\_}^2 $ & $2t_2t_3$
\\ $ \ [e^5] t_3  e^3  \beta_{8\_}^3 $ & $3t_2t_3$
\\ $ \ [e^5] t_5  e^5  \beta_{8\_}^5 $ &  \hspace{31pt}  $t_5$
\\ $[e^5]  (\beta_{8\_}-1 )$ & $5t_2t_3+t_5$
\\ \hline 
 $ \ [e^6] t_2  e^2  \beta_{8\_}^2 $ & $5t_2^3+2t_2t_4$
\\ $ \ [e^6] t_3  e^3  \beta_{8\_}^3 $ & \hspace{56pt}  $3t_3^2$
\\ \ $ [e^6] t_4  e^4  \beta_{8\_}^4 $ & \hspace{22pt} $4t_2t_4$
\\ $ \ [e^6] t_6  e^6  \beta_{8\_}^6 $ & \hspace{80pt} $t_6$
\\ $[e^6]  (\beta_{8\_}-1 )$ & $5t_2^3+6t_2t_4+3t_3^2+t_6$
\\
\hline
$ \ [e^7] t_2  e^2  \beta_{8\_}^2 $ & $12t_2^2t_3+2t_2t_5$
\\ $ \ [e^7] t_3  e^3  \beta_{8\_}^3 $ & \ \  $9t_2^2t_3$ \hspace{25pt} ${}+3t_3t_4$
\\ \ $ [e^7] t_4  e^4  \beta_{8\_}^4 $ & \hspace{68pt} $4t_3t_4$
\\ $ \ [e^7] t_5  e^5  \beta_{8\_}^5 $ & \hspace{37pt} $5t_2t_5$
\\ $ \ [e^7]t_7  e^7  \beta_{8\_}^7 $ & \hspace{102pt} $t_7$
\\ $[e^7]  (\beta_{8\_}-1 )$ & $21t_2^2t_3+7t_2t_5+7t_3t_4+t_7$
\\
\hline
$ \ [e^8] t_2  e^2  \beta_{8\_}^2 $ & $14t_2^4+14t_2^2t_4+\ 7t_2t_3^2+2t_2t_6$
\\ $ \ [e^8] t_3  e^3  \beta_{8\_}^3 $ & $\qquad \qquad \quad \quad  \ \ 21t_2t_3^2 \qquad \quad \  +3t_3t_5$
\\ \ $ [e^8] t_4  e^4  \beta_{8\_}^4 $ & $\qquad \quad 14t_2^2t_4 \qquad \qquad \qquad \qquad \qquad \ \ +4t_4^2$
\\ $ \ [e^8] t_5  e^5  \beta_{8\_}^5 $ & $\qquad \qquad \qquad \qquad \qquad \qquad \quad \quad 5t_3t_5$
\\ $ \ [e^8] t_6  e^6  \beta_{8\_}^6 $ & $\qquad \qquad \qquad \qquad \qquad \ \ 6t_2t_6$
\\ $ \ [e^8] t_8  e^8  \beta_{8\_}^8 $ & $\qquad \qquad \qquad \qquad \qquad \qquad \qquad \qquad \qquad \qquad t_8$
\\ $[e^8]  (\beta_{8\_}-1 )$ & $14t_2^4+28t_2^2t_4+28t_2t_3^2+8t_2t_6+8t_3t_5+4t_4^2+t_8$
\\
\hline
\end{longtable}
\end{center}
We omit rows with no terms.
Here, we've shown that $[e^k] (1 - \beta_{8\_} + \sum_{k \ge 2} t_k  e^k  \beta_{8\_}^k) = 0$ for natural $0 \le k <9$ so indeed $h(\beta_{8\_}) \mod e^9 = 0$.
\subsection{Face Layers}
For face layers, we set $v=e=1$ and we show that $h_3(\beta_{4\Delta 3}) \mod f^5 = 0$, with the subscripts $3$ indicating we are restricting ourselves to cubics, only $t_2$ and $t_3$ non-zero.
\begin{center}
\captionof{table}{Face Layering}\label{tab:face}  \addtocounter{table}{-1}
\begin{tabular}{ | p{.55cm} | p{3.4cm}| p{3.5cm} | p{3cm} |} 
\hline
 $f^F$  & $\beta_{4\Delta 3}-1 \qquad\qquad = {}$  & $t_2 f\beta_{4\Delta 3}^2  \qquad + {}$  & $ t_3 f \beta_{4\Delta 3}^3$
\\  \hline
$[f^0]$ & 0 & 0 & 0
\\ $[f^1] $&$ t_{2}+ t_{3} $ & $ t_{2}$ & $t_{3}$
\\ $[f^2] $&$  2  t_{2}^{2}  + 5  t_{2} t_{3}+ 3t_{3}^{2}   $&$   2  t_{2}^{2}  +2  t_{2} t_{3} $&$  3 t_{2} t_{3} + 3 t_{3}^{2}   $
\\ $[f^3]    $&  \mbox{ $5 t_{2}^{3} + 21 t_{2}^{2} t_{3} $}  \mbox{$\ \ {}+ 28 t_{2} t_{3}^{2} + 12 t_{3}^{3} $} 
&  $ 5 t_{2}^{3} +  12 t_{2}^{2} t_{3} +7 t_{2} t_{3}^{2}  $&$   9 t_{2}^{2} t_{3}+21  t_{2} t_{3}^{2}  + 12  t_{3}^{3}  $ 
\\ $[f^4] $ & \mbox{$14 t_{2}^{4} + 84 t_{2}^{3} t_{3} + 180 t_{2}^{2} t_{3}^{2} $}  \mbox{     $ {}+ 165 t_{2} t_{3}^{3} + 55 t_{3}^{4} $}
&$ 14 t_{2}^{4} + 56 t_{2}^{3} t_{3}
 + 72 t_{2}^{2} t_{3}^{2} + 30 t_{2} t_{3}^{3} $&$
28 t_{2}^{3} t_{3} + 108 t_{2}^{2} t_{3}^{2} + \ \  135 t_{2} t_{3}^{3} + 55 t_{3}^{4} $ \\
\hline
\end{tabular}
\end{center}
 Note $\Sb - 1$ at each face layer factors, in general with $t_2+t_3 + \ldots$ as a factor; the remaining factor is called the \textbf{Geode}~\cites{Wildberger2025, rubine2025hyper, rubine2025exer, Gessel2025, Zeilberger2025, Gossow2025}.

\section{Finite Layerings of Subdigons}  \label{sec:subdigons}

The polynomials we are examining are an associative and commutative projection of the algebra of multisets of subdigons.
A \textbf{subdigon} is a convex, planar, roofed polygon subdivided by non-crossing diagonals. 
The \textbf{roof} refers to a distinguished side of the polygon, drawn in red.
The inner, unsubdivided polygon with the roof as a side is called the \textbf{central polygon}.
The null subdigon, denoted $|$, is a subdigon with two vertices, one edge, no faces, and is of type $\m=[\ ]$, all $m_k=0$.
Its accounting monomial is $\psi(|)=1$.

\begin{figure}[H]
    \centering
    \includegraphics[width=0.6\linewidth]{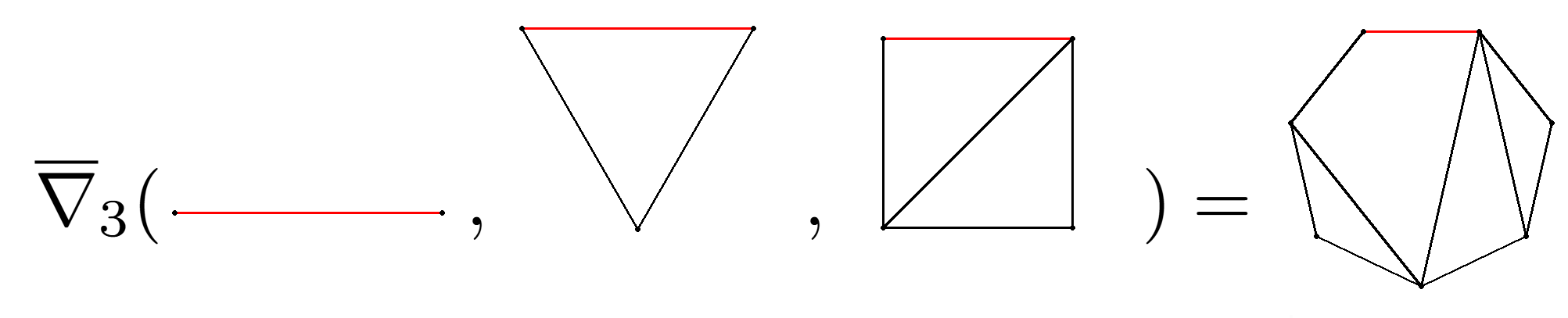}
    \caption{The $\tri_3$ operator creates a subdigon with a central quadrilateral}
    \label{fig:tri3}
\end{figure}

Subdigons have a recursive structure.
A subdigon is either the null subdigon $|$ or it is $\tri_k(s_1, s_2, \ldots, s_k)$ for some subdigons $s_1, \ldots, s_k$.
$\tri_k$ is a $k$-ary \textbf{paneling operator}; we have a family of them for natural $k \ge 2$.
$s = \tri_k(s_1, s_2, \ldots, s_k)$ is a new subdigon whose central polygon is a $k+1$-gon, with $k+1$ sides.
We designate a side as side zero, the new roof, then we proceed counterclockwise, designating side $1$ through side $k$. 
We adjoin $s_1$ by its roof to side $1$, $s_2$ by its roof to side $2$, etc., lastly adjoining $s_k$ by its roof to side $k$. 
For drawing, we normalize $s$ as a regular polygon, roof on top. That's permissible as we are only concerned with the combinatorial properties relating the various edges, faces and vertices of $s$ to its roof (figure \ref{fig:tri3}).

$\tri_k$ is non-associative and non-commutative; in fact for any non-null subdigon $s$ we can always recover the unique $s_1$, $s_2$, \ldots $s_k$ such that $s=\tri_k(s_1, s_2, \ldots, s_k)$.
Since $\tri_k$ adds a $k+1$-gon, and multiplication of the $\psi$s is addition of the type vectors, we get the basic but important identity:
\begin{equation} \label{eqn:psitrik}
\psi( \tri_k(s_1, s_2, \ldots, s_k) ) = t_k \psi(s_1)\psi(s_2)\cdots\psi(s_k)
\end{equation}
We extend the paneling operators to multisets of subdigons:
\begin{equation} \label{eqn:trikM}
\tri_k(M_1, \ldots, M_k) \equiv [ \ \tri(s_1, \ldots, s_k)\!: s_1 \in M_1, \ldots, s_k \in M_k \ ]    
\end{equation}
That extends the identity to multisets of subdigons per the linearity of $\Psi$:
\begin{equation} \label{eqn:PsitrikM}
\Psi( \tri_k(M_1, \ldots, M_k)) =  t_k \Psi(M_1) \Psi(M_2) \cdots \Psi(M_k)
\end{equation}
Let $\Sm_\m$ be the finite multiset of subdigons of type $\m$; we form the multiset $\Sm$ of all subdigons as layered by type:  
\begin{equation}
\Sm \equiv \sum_{\m \ge 0} \Sm_\m 
\end{equation}

\begin{equation}
\Psi(\Sm) =\sum_{\m \ge 0} \Psi(\Sm_\m) = \sum_{\m \ge 0} \mC \T^\m=\Sb  
\end{equation}

\begin{figure}
    \centering
    \includegraphics[width=.75\linewidth]{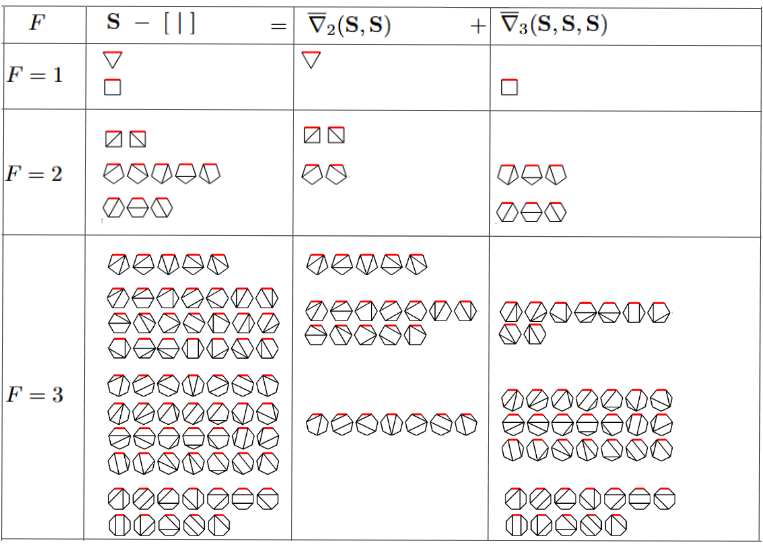}
    \caption{Face Layers of Cubic Subdigons}
    \label{fig:facesubd}
\end{figure}

We can write the subdigon grammar as a multiset equation:
\begin{equation} \label{eqn:Sm}
    \Sm = [ \ | \ ] + \tri_2(\Sm, \Sm) + \tri_3(\Sm, \Sm, \Sm) + \tri_4(\Sm, \Sm, \Sm, \Sm) + \ldots
\end{equation}
This says a subdigon is either the null subdigon, or it has a central triangle, or it has a central quadrilateral, etc.
That clearly applies at every layer.
For example, every subdigon with five edges has either a central triangle, quadrilateral or pentagon.  (One with a central hexagon will have at least six edges.)
It's of course true for vertex and face layers as well.
Figure \ref{fig:facesubd} illustrates the solution for face layers suitable for cubics (only subdividing into triangles and quadrilaterals).
Table \ref{tab:face} is essentially the $\Psi$ map applied to figure \ref{fig:facesubd}.

\section{Powers of Series Zeros} \label{sec:powers}

Applying $\Psi$ to both sides of equation \eqref{eqn:Sm}, we get 
\begin{equation}
\Sb = 1 + t_2 \Sb^2 + t_3 \Sb^3 + t_4 \Sb^4 + \ldots    
\end{equation}
i.e., $g(\Sb)=0$, $\Sb$ is a formal series zero of $g(\alpha) $ (equation \eqref{eqn:g}).

The question of powers of $\Sb$ arises naturally because each power appears as a term in the general geometric polynomial.

We saw from equation \eqref{eqn:Sm} that each such term comes from applying    $\Psi$ to multisets $\Ms_k = \tri_k(\Sm,\Sm, \ldots)$, subdigons with central $k+1$-gons.
So the generating series for subdigons with central $k+1$-gons must be:
\begin{equation} \label{eqn:Mk}
\Mb_k = \Psi(\Ms_k) = \Psi( \tri_k(\underbrace{\Sm,\Sm, \ldots, \Sm}_{k \textrm{ parameters}})) = t_k \Sb^k
\end{equation}
We see the powers of $\Sb$ are related to subdigons with particular central polygons.

In a similar context,
S. R. Mane~\cite{Mane2016} states, `it is conventional to solve for the $r$th powers of the roots,' and cites early examples from Lambert, 1770~\cite{Lambert1770}, Lagrange 1770~\cite{Lagrange1770}, and Fuss, 1795~\cite{Fuss1795}.

In the next section, we consider the simpler problem of powers of the Catalan generating series.
In section \ref{sec:raney1}, we describe Raney's~\cite{Raney1960} combinatorial derivation of the coefficients of $\Sb^r$.
Raney is working with strings describing functional compositions;
we give a simpler proof of Raney's essential lemma in section \ref{sec:raney2}.
That is followed by our derivation in terms of subdivided polygons in section \ref{sec:central}.

\section{A Recurrence for Powers of the Catalan Generating Series} \label{sec:cat}

We start simply, with powers of $\Tb$, where $\Tb$ is the generating series for the Catalan numbers:
\begin{equation}
    \Tb = \sum_{n \ge 0} C_n t^n
\end{equation}
$\Tb = \Sb[t, 0, 0, \ldots]$ satisfies the quadratic equation:
\begin{equation}
    \Tb = 1 + t \Tb^2
\end{equation} 
We can view the quadratic equation as a rewrite rule that allows us to reduce the power of $\Tb$.
\begin{theorem}[Powers of the Catalan generating series]
    $$C^{(r)}_m=\frac{r}{2m+r}\binom{2m+r}{m}$$
where $C^{(r)}_m=[t^m]\Tb^r$, the coefficient of $t^m$ in $\Tb^r$.
\end{theorem}

\begin{proof}
For the base case, $r=1$, clearly $\displaystyle C_m^{(1)} = \dfrac{1}{2m+1}\binom{2m+1}{m}=\dfrac{1}{m+1}\binom{2m}{m}=C_m$, and of course
that is correct because $[t^m]\Tb^1 = [t^m]\Tb = C_m$.

Now we assume that for every $r < R$, $C^{(r)}_m=\frac{r}{2m+r}\binom{2m+r}{m}$. Then, we need to prove that $C^{(R)}_m=\frac{R}{2m+R}\binom{2m+R}{m}$. 
To do so, we develop a recurrence for powers of $\Tb$, from the defining quadratic
$ \Tb = 1 + t \Tb^2$.
Our goal is to find $P_r=P_r(t)$ and $Q_r=Q_r(t)$, polynomials in $t$, that let us express:
\begin{equation} \label{eqn:T}
t^{r-1} \Tb ^r=P_r \Tb +Q_r 
\end{equation}
For $r=1$, $t^0\Tb=\Tb$ so $P_1(t)=1$ and $Q_1(t)=0$.
For $r=2$, $t \Tb ^2=\Tb-1$, so $P_2(t)=1$ and $Q_2(t)=-1$. 
Multiplying equation \eqref{eqn:T} by $t\Tb$,
\begin{equation}
t^r\Tb^{r+1}=P_rt\Tb^2+Q_rt\Tb
=P_r(\Tb-1)+Q_rt\Tb
=(P_r+tQ_r)\Tb-P_r
\end{equation}

Thus, we have that
$Q_{r+1}=-P_r$. 
\begin{align}
t^r\Tb^{r+1}&=(P_r- t P_{r-1})\Tb-P_r
\end{align}

We found a recurrence for $P_r(t)$: 
\begin{align} \label{eqn:reccurence}
P_0(t)=0, \ P_1(t)=1, 
\ P_r(t) & =P_{r-1}(t)-t P_{r-2}(t) 
\end{align}

We conclude the degree in $t$ of $P_r(t)$ is $\lfloor \frac{r-1}{2} \rfloor$.
We introduce $m$ with the property $1 \le r-1 \le m$.
We note $\deg P_{r-1}(t) < r-1$ so $[t^m]P_{r-1}(t)=0$.
\begin{align}  \label{eqn:CatalanPowertoP}
    C^{(r)}_{m-r+1} &=[t^{m-r+1}] \Tb^r = [t^m]t^{r-1}\Tb^r
=  [t^m] P_r \Tb +  [t^m] Q_r 
\\& \nonumber  = [t^m] P_r \sum_{k \ge 0} C_k t^k - [t^m] P_{r-1}
= \sum_{0\le k\le m}{C_{m-k} \,[t^k] P_r(t)}  
\end{align}

Thus, we get a formula for the Catalan Powers in terms of $P_r(t)$. We can proceed by applying our recurrence equation \eqref{eqn:reccurence} to the $P_r(t)$ in the summation above.

\begin{align} 
 C^{(r)}_{m-r+1}  &= \sum_{0\le k\le m}{C_{m-k} \, [t^k] (P_{r-1}(t)-tP_{r-2}(t))}
 \\  \nonumber & = \sum_{0\le k\le m}{C_{m-k} \, [t^k] P_{r-1}(t)}-\sum_{0\le k\le m}{C_{m-k} \, [t^{k-1}]P_{r-2}(t)}
 \\  \nonumber & = \sum_{0\le k\le m}{C_{m-k} \, [t^k] P_{r-1}(t)}-\sum_{0\le k\le m-1}{C_{m-1-k} \, [t^{k}]P_{r-2}(t)}
 \\  C^{(r)}_{m-r+1}  &=  C^{(r-1)}_{m-r+2}-C^{(r-2)}_{m-r+2}
\end{align}
where we get the last expression from equation \eqref{eqn:CatalanPowertoP}. We set  $m'=m-r+1$  and drop the primes to get:
\begin{equation}
    C^{(r)}_m=C^{(r-1)}_{m+1}-C^{(r-2)}_{m+1}
\end{equation}
Setting $r=R$ and applying the induction assumptions  
$C^{(r)}_m=\frac{r}{2m+r}\binom{2m+r}{m} = \frac{r(2m+r - 1)!}{ (m+r)! m! }$ for $r < R$,
\begin{align}
C^{(R)}_{m} &= C^{(R-1)}_{m+1}-C^{(R-2)}_{m+1}
 \\
&=\frac{(R-1)(2m+R)!}{(m+R)!(m+1)!}-\frac{(R-2)(2m+R-1)!}{(m+R-1)!(m+1)!}
\nonumber \\ C^{(R)}_{m}  &= \frac{R(2m+R-1)!}{(m+R)!m!}
= \frac{R}{2m+R}\binom{2m+R}{m}
\end{align}
That is the induction step required for our proof.
\end{proof}

\section{Raney's Derivation of Hyper-Catalan Power Numbers} \label{sec:raney1}

That algebraic proof was a bit complicated, and the generalization to $\Sb^n$ looks even more complicated.
We now review a combinatorial derivation due to Raney, 1960~\cite{Raney1960}, also a bit complicated.

Raney derives the formula for the hyper-Catalan numbers and the hyper-Catalan power numbers (coefficients of $\Sb^n$) at the same time, by comparing lists of functional composition expression strings (in a Polish notation with an obvious bijection to plane trees) to a general string permutation of a type given by a multinomial coefficient.

In the realm of subdivided polygons, 2-gons are somewhat problematic, and not necessary for their solution, so $m_1$ and $t_1$ were omitted by Wildberger and Rubine.
Raney allows unary functions, which map to unary nodes in plane trees,
making Raney's array more general than the hyper-Catalans.
His expressions generalize Wildberger and Rubine's to include $m_1$ and $t_1$ in the obvious way; 

 see also Gessel, 2025~\cite{Gessel2025}.

Raney defines several terms as follows.
A \textbf{string} $\sigma$ is a finite sequence of natural numbers $a_1a_2a_3 \cdots a_m$.
The \textbf{length} of this string is $m$.
The \textbf{rank} of this string is $r(\sigma)=\sum_{k=1}^{m}(a_k-1)$. 
A string is a \textbf{word} if it is a single $0$ or the string $n \alpha_1 \alpha_2 \cdots \alpha_n$ for a natural number $n>0$, where each string $\alpha_i$ is a word.

Raney relies on a theorem of Rosenbloom~\cite{Rosenbloom1950} that the string   $a_1a_2a_3 \cdots a_m$ is a word
precisely when it has rank $-1$ and no \textbf{proper prefix} $a_1a_2\cdots a_l$ for $l<m$ has negative rank.
He further relies on an inductive extension of Rosenbloom's theorem to lists of $n$ words  (given without proof): a string is a list of $n$ words precisely when the rank is $-n$ and no proper prefix has rank less than or equal to $ -n$.

A word is Raney's specification of a well-formed functional composition, a string with a tree-like recursive structure.
As an example,
representing each natural number as a digit,
the functional composition $f(0,f(0,f(0,f(0),0)))$ corresponds to the word $202030100$.

\subsection{Raney's Lemma}

Raney begins with a lemma (his Theorem 2.1):
for natural $n>0$, each string $\sigma$ of rank $-n$ has precisely $n$ cyclical rotations that are lists of $n$ words. 
We sketch his proof, the heart of his argument.

He proceeds by considering each cyclic rotation of a string $\sigma$ of rank $-n$, and among all those he finds the prefix $\alpha$ whose rank is maximal, and whose length is maximal among those prefixes with maximal rank. He calls $\sigma_1$ the rotation that gives rise to $\alpha$ as maximal prefix. Considering every proper prefix $\tau$ of $\sigma_1$, there are two cases: $\alpha=\tau \gamma$ or $\sigma_1=\alpha \delta \epsilon$ where $\alpha\delta=\tau$ and $\epsilon$ is not empty.  

Raney wants to show the rank of $\tau$ is bigger than $-n$.
In the first case, $\alpha=\tau \gamma$ so 
$r(\alpha)=r(\tau)+r(\gamma)$.
Raney argues: Assuming $r(\tau) \le -n$ then $r(\alpha) \le -n+r(\gamma)<r(\gamma)$, a contradiction, so $r(\tau)>-n$ in this case.  
We're not sure why he couldn't just skip the indirect argument and write
 $r(\tau)=r(\alpha)-r(\gamma) \ge 0 > -n$ as $\alpha$ has maximal rank and $\gamma$ is a prefix of some rotation.

In the second case,  $\sigma_1=\alpha \delta \epsilon$ where $\alpha\delta=\tau$ and $\epsilon$ is not empty. 
Again Raney assumes $r(\tau) \le -n$.
Here $-n = r(\sigma_1)=r(\tau \epsilon) = r(\tau)+r(\epsilon) \le -n + r(\epsilon)$. Then $r(\epsilon) \ge 0$ so $r(\alpha) \le r(\epsilon) + r(\alpha)= r(\epsilon \alpha)$, a contradiction because $\alpha$ is maximal rank and length, so $\epsilon\alpha$ cannot have a larger rank or the same rank and be longer.

So $\sigma_1$ has rank $-n$ and no proper prefix has rank as small.
From the generalization of Rosenbloom's theorem, that would be sufficient to
conclude $\sigma_1$ is a list of $n$ words.
But Raney continues the argument.

Raney shows $\sigma_1$ factors as $\sigma_1 = \alpha_1 \alpha_2 \cdots \alpha_n$ where $\alpha_1$ is the first prefix of $\sigma_1$ with rank $-1$ so $\sigma_1 =\alpha_1 \beta_1$; $\alpha_2$ is similarly the first prefix of $\beta_1=\alpha_2 \beta_2$ of rank $-1$, etc. $\beta_n$ is empty because $\sigma_1$ has no proper prefixes of rank $-n$.
Therefore $\sigma_1$ is a list of $n$ words, as are all $n$ rotations beginning with an $\alpha_i$.
(This is applying Rosenbloom's theorem to each $\alpha_i$.)

To show there are no additional lists of $n$ words, he imagines splitting one cyclic rotation $\sigma_2$ in the middle of a word:  $\sigma_2 =\beta \gamma \delta$, where $\beta$ and $\delta$ are non-empty and $\delta\beta=\alpha_i$ for some $i$.  Then $r(\delta) \ge 0$ because it is a prefix of a word, so $r(\beta\gamma) \le -n$, so $\sigma_2$ cannot be a list of $n$ words because its proper prefix $\beta\gamma$ violates the rank condition.
(This appears to be an application of the generalization of Rosenbloom's theorem.)
So Raney's lemma is proven.

\subsection{Raney counts \textit{n} word lists of a type} \label{sec:Raney22}

Raney defines $L(n, m_1, m_2, m_3, \ldots)$ to be the number of word lists with $n$ words, and a total of $m_1$ $1$s, $m_2$ $2$s, $m_3$ $3$s, etc.
Once these are specified, the number of zeros is determined by the condition that the rank be $-n$ as $m_0=n + m_2 + 2m_3 + \ldots$.
The length of each such word list is $m = m_0+m_1 +m_2+ \ldots= n + m_1 + 2m_2 + 3m_3 + \ldots $.
The number of strings of that length with the given composition (each of rank $-n$) is given by the multinomial coefficient
$ \binom{m}{m_0 , m_1,  m_2, m_3, \ldots} $.
Raney's lemma tells us exactly $n$ cyclic rotations of each such string are lists of $n$ words.
There are $m$ possible cyclic rotations, so 
Raney concludes (his Theorem 2.2):
\begin{equation} \label{eqn:L}
m L(n, m_1, m_2, m_3, \ldots) = n \binom{m}{m_0, m_1,  m_2, m_3, \ldots}
\end{equation} 

so his generalization of the hyper-Catalan power coefficient (including $m_1$) is:
\begin{align}
L(n, & m_1, m_2, m_3, \ldots)
\\ \nonumber &= \dfrac{n}{n + m_1 + 2m_2 + 3m_3 + \ldots} \binom{n + m_1 + 2m_2 + 3m_3 + \ldots}{n + m_2 + 2m_3 + \ldots  , m_1,  m_2, m_3, \ldots}
\end{align} 

Clearly $n=1$, a list with a single word, gives (the generalization of) the hyper-Catalan numbers $\mC=L(1,m_1,m_2,m_3,\ldots)$.
We consider raising its generating function $\Sb=\sum_{\m \ge 0} \mC \T^{\m}$ to the $n$th power.
(Here we conflate $\Sb$ with this more general series that includes a $t_1$.)
Raney observes that the coefficient $L(n, m_1, m_2, \ldots)$, the number of lists of $n$ words with the indicated composition, corresponds to the ways to collect like terms of $\T^\m$ in $\Sb^n$; all the ways the type vectors of the $n$ factors add up to $\m$.
\begin{equation}
[t_1^{m_1}t_2^{m_2}t_3^{m_3}\cdots] \Sb^n = L(n, m_1, m_2, m_3, \ldots) 
\end{equation}

\subsection{Examples}

Consider type $[m_1, m_2, m_3] = [0, 2, 1]$ and power $n=1$.
We have $m_0=1+2(1)+1(2)=5$ 0's, $2$ 2's, and $1$ 3's in our list of 1 word.
Our word length is $m=1+0(1)+2(2)+1(3)=8$ so
there are $\binom{8}{5, 0, 2, 1}=168$ different permutations that can be formed. Out of those, exactly $\frac{1}{8}$ of them (21) are valid words, those being:

\noindent
20203000,
20230000,
20300200,
20302000,
20320000,
22003000,
22030000,
\newline
22300000,
23000200,
23002000,
23020000,
23200000,
30020200,
30022000,
\newline
30200200,
30202000,
30220000,
32000200,
32002000,
32020000,
32200000.

Notice that no two words are cyclic permutations of one another.
Each of these corresponds to a subdigon in figure \ref{fig:C21}; can you match them?

For a second example, consider $n=3, m_1=1, m_2=1$. So then $m_0=3+1(1)=4$. Our word length is $m=3+1(1)+1(2)=6$ so
there are $\binom{6}{4, 1, 1}=30$ different permutations that can be formed. Out of those, exactly $\frac{3}{6}$ of them (15) are valid lists of three words, those being:

\noindent
001200,
002010,
002100,
010200,
012000,
020010,
020100,
021000,
100200,
\newline
102000,
120000,
200010,
200100,
201000,
210000.

We can verify this.
$$[t_1^1t_2^1] (\Sb[t_1, t_2, \ldots])^3= [t_1^1t_2^1](1 + t_1 + t_2  + 3t_1 t_2 + \ldots)^3 = 15$$

\section{Clarifying Raney} \label{sec:raney2}

We can prove Raney's lemma in a simpler, more intuitive way.
By introducing the composition of the string at the outset, we no longer require Rosenbloom's Theorem or Raney's generalization. 

\begin{theorem}[Raney's lemma]
For natural $n>0$, each string of rank $-n$ has precisely $n$ cyclical rotations that are lists of $n$ words.   
\end{theorem}
\begin{proof}
We have a rank $-n$ string $\sigma$ which we arrange in a circle.
It has $m_0$ 0's, $m_1$ 1's, $m_2$ 2's, $m_3$ 3's, etc.
Expressing the rank in terms of the type, that means $-n=\sum_{k \ge 0} (k-1) m_k $.
So we must have $m_0=n+m_2+2m_3+\cdots$ 0's, each of which adds $-1$ to the rank.

We describe an algorithm to identify the $n$ words formed by circular $\sigma$. Each $0$ is an \textbf{identified word}; there's no need to mark it further. We start by finding all strings of the form $i$ followed by $i$ zeroes for $i>0$; each of these is an identified word; we mark them by enclosing each in parentheses.
Note there is no ambiguity, no possibility of overlap.
We repeat this process, at each step grouping any instance of the number $i$ followed by $i$ identified words (including $0$'s) to form a bigger identified word. 

We claim that once this process terminates, we are left with $n$ identified words.

To show this, we extend our notion of type to strings with identified words by including the count of identified words in $m_0$, and not counting their elements in the other $m_i$. 
When we perform the above process of identifying words, we group together $1$ of the $m_i$'s and $i$ of the $m_0$'s to form a new identified word, counted in $m_0$. So we decrease $m_i$ by $1$ and $m_0$ by $i-1$. When we do this, $\sum_{k \ge 0}{(k-1)m_k}$ doesn't change.
Thus the rank is invariant, remaining $-n$.
Accounting-wise, it's as if we've replaced each identified word (always rank $-1$) with a zero (also  rank $-1$).

Say, for the sake of contradiction, that the process terminated before we reached a list of identified words, i.e. terminated with $m_i>0$ for at least one $i \ge 1$. Then every $1$ must have less than $1$ identified word immediately after it, every $2$ must have less than $2$ identified words immediately after it, and so on. Each substring consisting of $i$ followed by less than $i$ identified words would have a nonnegative rank. Thus the sum, the rank of the entire string $\sigma$, would be nonnegative, a contradiction as $\sigma$'s invariant rank is $-n$, $n>0$.

Thus, this process only terminates once we deplete all $m_i$ for $i \ge 1$. From our invariant, we know that $\sum_{k \ge 0 }{(k-1)m_k}=-n$, so at the end we must have $m_0=n$, a list of exactly $n$ words.

This algorithm finds all substrings that are words, and deterministically leaves us with a list of $n$ words.
We conclude that there are exactly $n$ cyclic rotations of a string with rank $-n$ that are a list of $n$ words.
\end{proof}

Raney's counting formula (equation \ref{eqn:L}) follows from his lemma, as already shown in section \ref{sec:Raney22}.

\subsection{Example}
Each a step in the algorithm gives a cyclic string, arranged in a line for presentation, so a parenthesized identified word may wrap around the string.

In our example, $m_0=12, m_1=3, m_2=1, m_3=2, m_4=1$, so rank $-n=-4$.

0030130010001000420

0)0301300(10)00(10)004(20

0)0301(300(10))00(10)004(20

0)030(1(300(10)))00(10)004(20

0)0(30(1(300(10)0)))0(10)004(20

0)0(30(1(300(10)0)))0)(10)00(4(20

There are no more moves so the process terminates and we are left with a list of $4$ identified words: (10), 0, 0, (4(200)0(30(1(300(10)0)))0).

\section{Hyper-Catalan Power Numbers via Central Polygons} \label{sec:central}
$\mC$ counts the number of subdigons of type $\m=[m_2,m_3,m_4, \ldots]$ and has generating series
$\Sb=\sum_{\m\ge \mathbf{0}} \mC \T^{\m}$.
In this section we're after the closed form of
$\mC^{(r)} \equiv[\T^\m]\Sb^r$, which we're calling the \textbf{hyper-Catalan Power  Numbers}.
We can get it from Raney's result in section \ref{sec:raney1} as $L(r, 0, m_2, m_3, \ldots)$ but we develop our own derivation here.

In section \ref{sec:powers} we defined  $\Ms_k = \tri_k(\Sm,\Sm, \ldots)$ to be the multiset of subdigons with central $k+1$-gons, with generating series $\Mb_k=\Psi(\Ms_k) = t_k \Sb^k$.  We want to count these by type.
\begin{figure}
    \centering
    \includegraphics[width=1\linewidth]{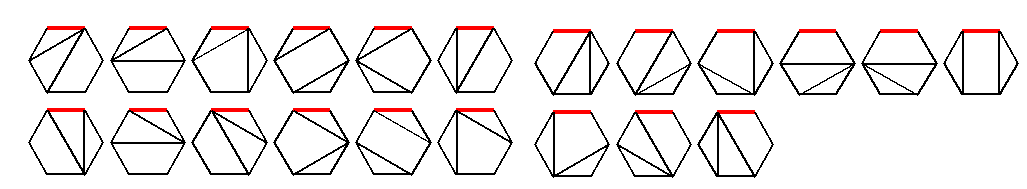}
    \caption{$C[2,1]=21, V=2+1(2)+2(1)=6, E=1+2(2)+3(1)=8, F=2+1=3$}
    \label{fig:C21}
\end{figure}
We seek $C_{\m,r}$, the number of subdigons of type $\m$ with a central $r+1$-gon.  
In the example in figure \ref{fig:C21} we see $C[2,1]=21$ subdigons with two tris and a quad.  12 of them have central tris, 9 have central quads.  We have $\mC = \dfrac{(\mE-1)!}{(\mV-1)! \, \m!}$.  \
We note numerator factor $\mE-1 = 2m_2+3m_3+\ldots$ gives us a contribution of $2(2)=4$ from tris and $3(1)=3$ from quads, and indeed $\frac{4}{4+3}\mC=12, \frac{3}{4+3}\mC=9$.
It's easy to guess:

\begin{theorem} The number of subdigons of type $\m$ with a central $r+1$-gon is:
$$C_{\m,r} = \dfrac{r \, m_r \, \mC}{\mE -1}$$
\end{theorem}
\begin{proof}
Since $\mE= 1 + \sum_{k \ge 2} k \, m_k$ (equation \eqref{eqn:CVMF}), 
each $r+1$-gon contributes $r$ new edges to a subdigon.

Among all the subdigons of type $\m$, each edge is equally likely to be the roof.
So $\dfrac{r m_r}{\sum_k k m_k} = \dfrac{r m_r}{\mE-1}$ of the $\mC$ subdigons have a central $r+1$-gon.
\end{proof}

We define basis vectors $\vec{j} \equiv [0,0,\ldots,0,1]$ where there are $j-2$ zeros.
From equation \eqref{eqn:Mk} we know the generating function for subdigons of a type with a central $r+1$-gon is:
\begin{align}
\Mb_r & =\Psi(\Ms_r) = t_r \Sb^r = \sum_{\n \ge 0} C_{\n,r} \T^{\n} 
\\\Sb^r &= t_r^{-1} \sum_{\n \ge 0} C_{\n,r} \T^{\n} = \sum_{\n \ge 0} C_{\n,r} \T^{\n-\vec{r}}
\\ \mC^{(r)}  &= [\T^\m] \Sb^r   =  C_{\m+\vec{r},r}
\end{align}
That's the proof of the following theorem.
\begin{theorem} (Hyper-Catalan power coefficients via central polygons) The hyper-Catalan power coefficient $\mC^{(r)} \equiv [\T^\m]\Sb^r$ is given by the number of subdigons of type $\m+\vec{r}$ with a central $r+1$-gon:
   $$ \mC^{(r)} = C_{\m + \vec{r},r}$$
\end{theorem}
It's analogous to Raney, where a subdigon with a central $r+1$-gon is playing the role of a list of $r$ subdigons, the ordered children of the central $r+1$-gon.

\begin{theorem}(Closed form of the hyper-Catalan power coefficients)
$$\mC^{(r)} \equiv[\T^\m]\Sb^r = \dfrac{r ( r-2 + \mE  )!}{(r-2  +\mV)! \, \m!}  
$$
\end{theorem}
\begin{proof}
Just some algebra to work out $ C_{\m + \vec{r}, r}$. 
\begin{align}
\mC^{(r)} &= C_{\m + \vec{r}, r}
= \dfrac{r (m_r+1) C_{ \m + \vec{r}} }{E_{\m + \vec{r}} -1} 
\\ \nonumber  &
= \dfrac{r  (m_r+1) }{E_{\m + \vec{r}} \, -1} \cdot \dfrac{(  E_{\m + \vec{r}} \, -1 )! }{(  V_{\m + \vec{r}} \, -1 )!  (\m+\vec{r})!  }
\\ \nonumber  &= \dfrac{r  (m_r+1) ( \mE + r -2 )! }{(  \mV + r-2)!  \m! (m_r+1)  }
= \dfrac{r  ( r -2+ \mE  )! }{(  r-2 + \mV )!  \m! }
\end{align}
\end{proof}
We note the expanded form:
\begin{equation}
    \mC^{(r)}  = \dfrac{r ( r-1 + 2m_2 + 3m_3 + 4m_4 + \ldots )!}{(r + m_2 + 2m_3 + 3m_4 + \ldots)! \,  m_2! \, m_3! \,  m_4! \cdots}
\end{equation}

\section{Conclusion}
The hyper-Catalan generating series $\Sb[t_2, t_3, t_4, \ldots]$, a zero of $g(\alpha)$, appears aggressively infinite, not only in the number of terms in the sum, but also in the unbounded number of independent variables $t_2, t_3, \ldots \ $.
There is a term for every possible product of natural powers of the $t_k$.
Constraining $g(\alpha)$ to be a polynomial at least restricts us to a finite number of $t_k$, but each associated $m_k$ still ranges from 0 to infinity.
By layering the series by level, where a level is a maximum number of vertices, number of faces, or number of edges, we get a truncated version of $\Sb$ that is still a zero of $g(\alpha)$, but only when we discard all terms of higher level than we are considering.
For face layering we need to further bound the degree of the equation for finite layers.
In each of the three layerings, we are able to interpret the series solution as a family of finite identities.

In addition to the solution $\Sb$, methods such as Lagrange Inversion~\cites{Lagrange1770, Gessel2016} often give expressions for the coefficients of the powers $\Sb^r$.
Raney offers a combinatorial explanation for the coefficients of $\Sb^r$, as we recount.
Unfortunately Raney is a bit difficult, even in our recounting, so we have also produced our own proof of Raney's key lemma, as well as a derivation involving subdigons with central polygons, which we hope are easier to understand. 
A subdigon with a central $r+1$-gon is isomorphic to Raney's list of $r$ words, but easier to count.

Raney, 1960~\cite{Raney1960} counted words and lists of words of a type; his words are strings of natural numbers explicitly isomorphic to rooted plane trees.
Rhoades, 2011~\cite{Rhoades2011} lists more (what we call) hyper-Catalan objects, including  typed noncrossing set partitions, typed nonnesting set partitions, typed Dyck paths, and rooted plane trees of a type (first counted by Tutte in 1964~\cite{Tutte1964}).
Shuetz and Wheildon, 2016~\cite{Schuetz2016} add typed subdivided polygons (what we've called subdigons), originally counted by Erd\'elyi and Etherington in 1940~\cite{Erdelyi1940}.
Kreweras, 1972~\cite{Kreweras1972} counts non-crossing partitions of a cyclic graph (similar to Chu, 1987~\cite{CHU198791}), and
Gessel, 2025~\cite{Gessel2025} shows certain lattice path types are counted by the hyper-Catalans (and others by the Geode!).
Each bijection affords different proof techniques; Raney's neat cyclic rectangle of words admit different analyses than do subdigons or rooted plane trees or lattice paths.
Furthermore, each of the hundreds of Catalan number  bijections~\cite{Stanley2015} may have an associated hyper-Catalan generalization to be considered.

Wildberger and Rubine's series solution, and our finite interpretation of it, offer a different understanding of what a zero of a polynomial really is.
The arithmetic takes place in a larger realm, where the multivariate polynomials and power series are themselves the arithmetic objects, the generalizations of numbers. 
In this context, a polynomial zero is not just a complex number, especially not just one that requires an infinite number of steps to compute.
It's a power series, that, even though is an ongoing expression, has a finite interpretation, as we have shown.

\bibliographystyle{vancouver}
\bibliography{HyperCatBib.bib}

\end{document}